\documentclass{article}     
                         
\usepackage{amssymb}

\usepackage{enumerate}
\usepackage{amsthm,amsmath,rotate,enumerate,mathrsfs,amsfonts,latexsym,amssymb,cancel,accents,graphicx,algorithm2e,csquotes,turnstile,url,stackengine,float,tikz,hyperref}
\usepackage{scalerel}
\usepackage{setspace}

\usetikzlibrary{backgrounds}

\def\newop#1{\expandafter\def\csname #1\endcsname%
     {\mathop{\rm #1} \nolimits}}

\newcommand\nprove{\mathbin{\cancel{\vdash}}}

\makeatletter
\newcommand*\bigcdot{\mathpalette\bigcdot@{.5}}
\newcommand*\bigcdot@[2]{\mathbin{\vcenter{\hbox{\scalebox{#2}{$\m@th#1\bullet$}}}}}
\makeatother

\newtheorem{Theorem}{Theorem}

\newtheorem{Lemma}{Lemma}

\newtheorem{Proposition}{Proposition}
\newtheorem{Fact}{Fact}

\theoremstyle{definition}
\newtheorem{Definition}{Definition}

\newtheorem*{Main}{Main Theorem} 

\def\eq{\leftrightarrow}
\def\imp{\rightarrow}
\def\ex{\exists}
\def\all{\forall}

\def\Pr{\textup{\texttt{Pr}}}
\def\PA{{\mathsf{PA}}}

\def\Prf{\textup{\texttt{Prf}}}
\def\Formula{\textup{\texttt{Formula}}}

\def\Proof{\textup{\texttt{Proof}}}
\def\Con{\textup{\texttt{Con}}}

\def\s{\textup{\texttt{s}}}
\def\sub{\textup{\texttt{sub}}}

\def\impd{\overset{\bigcdot}{\rightarrow}}
\def\eqd{\overset{\bigcdot}{\leftrightarrow}}
\def\negd{\overset{\bigcdot}{\neg}}
\def\lord{\overset{\bigcdot}{\lor}}

\def\landd{\overset{\bigcdot}{\land}}

\def\ul{\ulcorner}
\def\ur{\urcorner}

\def\Oca{\omega\text{-}\Con}

\def\l{\langle}
\def\r{\rangle}
\def\EA{\mathsf{EA}}
\def\numeral{\texttt{numeral}}
\def\exp{\mathsf{exp}}

\def\varr{\textup{\texttt{-var}}}

\def\Tr{\textup{\texttt{Tr}}}

\begin{document} 
\title{$\omega$-consistency for Different Arrays of Quantifiers}


 




\author{Paulo Guilherme Santos\footnote{CMAFcIO, UL; ISCAL, IPL. Email: \href{mailto:paulo.g.santos.logic@gmail.com}{paulo.g.santos.logic@gmail.com} }}
\date{}
\maketitle

\begin{abstract} We study the formalized $\omega$-consistency statement by allowing the occurrence of different arrays of quantifiers in it. In more detail, for an array of quantifiers $\vec{Q}$, we introduce the $\omega$-consistency statement $\Oca^{\vec{Q}}_T$. In our framework, $\Oca_T$, the formalized $\omega$-consistency, is simply $\Oca^{\all}_T$. We prove that for some specific arrays of quantifiers we get consistency statements that are $S$-equivalent to the original $\omega$-consistency statement ($S$ denotes the basis theory to develop metamathematics), namely
\begin{align*} S\vdash\Oca_T\eq \Oca^{\all^{n+1}}_T\eq\Oca^{\all^{n+1}\ex^m}_T\eq\Oca^{\ex\all\ex^n}_T.
\end{align*}
We end our paper by creating a theory of truth that proves each $\Oca^{\vec{Q}}_T$-statement.
\end{abstract}

\section{Introduction}  
$\omega$-consistency was first introduced by K. G\"odel in \cite{godel1931formal} as a useful condition to prove his First Incompleteness Theorem. We recall that a theory $T$ is \emph{$\omega$-inconsistent} if there is a formula $\varphi(x)$ such that
\begin{description}
\item[\textbf{$\omega$-inc.1:}] For all $n\in \mathbb{N}$, $T\vdash \varphi(\overline{n})$, and
\item[\textbf{$\omega$-inc.2:}] $T\vdash \ex x.\neg \varphi(x)$;
\end{description}
and \emph{$\omega$-consistent} otherwise (see \cite{boolos1995logic} and \cite[p.\hspace{0.1cm}851]{Barwise1993} for further details on $\omega$-consistency). Inside a theory of arithmetic, this form of consistency can be represented by
\begin{align*}\all 1\varr(y). (\all x.\Pr_T(\s(y,x)))\imp\neg \Pr_T(\overbrace{\ex x_0.}^{\bigcdot}\negd \sub(y,\ul x_0\ur)),
\end{align*}
where, as usual, $\Pr_T$ denotes provability in $T$, $\s(y,x)$ denotes the ``the substitution in the formula (coded by) $y$ by the numeral of $x$'', $\sub(y,\ul x_0\ur)$ represents ``the substitution of the first variable in the formula (coded by) $y$ by the variable $x_0$'', $\overbrace{\ex x_0.}^{\bigcdot}$ represents (computable) function that, given as input the code of a formula, outputs the code of the formula obtained from the first one by placing `$\ex x_0.$' as prefix\footnote{We use distinctly the quotations ``$\cdots$'' for idiomatic expressions, and the quotations `$\cdots$' for arrays of text.}, and $n\varr(y)$ expresses ``the formula (coded by) $y$ has exactly $n$ free variable''. By defining $\overline{\ex}:=\all$ and $\overline{\all}:=\ex$ we can express $\omega$-consistency as
\begin{align*} \all 1\varr(y).(\all x.\Pr_T(\s(y,x)))\imp\neg \Pr_T(\overbrace{\overline{\all} x_0.}^{\bigcdot}\negd \sub(y,\ul x_0\ur)),
\end{align*}

In our paper, we study forms of $\omega$-consistency where we allow the occurrence of different arrays of quantifiers besides the one that we have just presented. More specifically, for an array of quantifiers $\vec{Q}$, we define
\begin{align*} \Oca^{\vec{Q}}_T:=\all |\vec{Q}|\varr(y). (\vec{Q}\vec{x}.\Pr_T(\s(y,\vec{x})))\imp\neg \Pr_T(\overbrace{\overline{\vec{Q}}\vec{x_0}.}^{\bigcdot}\negd \sub(y,\ul\vec{x_0}\ur)),
\end{align*}
where $|\vec{Q}|$ denotes the number of quantifiers in the array $\vec{Q}$. Clearly, $\Oca^{\vec{Q}}_T$ generalizes the standard $\omega$-consistency; in fact, $\omega$-consistency corresponds to $\Oca^{\all}_T$. As our main result, we show that there are different array of quantifiers $\vec{Q}$ (besides the obvious one, namely $\all$) that give rise to statements that are provably equivalent to $\omega$-consistency; in other words, $\omega$-consistency can be conceived using different quantifiers. More specifically, we prove
\begin{align*} S\vdash\Oca_T\eq \Oca^{\all^{n+1}}_T\eq\Oca^{\all^{n+1}\ex^m}_T\eq\Oca^{\ex\all\ex^n}_T,
\end{align*}
where $S$ denotes our weak basis theory used to perform metamathematics, and $Q^n$ denotes $\overbrace{Q\cdots Q}^{\text{$n$ times}}$. We also develop a theory of truth that proves each and every $\Oca^{\vec{Q}}_T$-statement. 

We end our introduction to the topic of this paper by giving a concrete example of our main result. It entails, for example, that the usual $\omega$-consistency could be stated in the following way:
\begin{quotation} For every two variable formula $\varphi(x,y)$, if for all $n\in \mathbb{N}$ there is $m\in \mathbb{N}$ such that $T\vdash \varphi(\overline{n},\overline{m})$, then $T\nprove \ex x. \all y. \neg \varphi(x,y)$.
\end{quotation}

\section{Preliminaries}  
Throughout this paper, $S$ and $T$ stand for consistent\footnote{Throughout our paper, whenever we mention the fact that a formula is true, we of course mean true in $\mathbb{N}$.} theories that include $\EA$ (this theory is $\mathsf{I}\Delta_0+\exp$, where $\exp$ denotes the totality of exponentiation; see \cite{BEKLEMISHEV1997193} for details on this theory); we assume $S\subseteq T$, $S$ is a potentially weak theory, and $T$ a potentially strong theory (for example $\PA\subseteq T$); we assume that $S$ has a \emph{modicum} of arithmetic required to state our results, we emphasize the requirements as we go: usually, the \emph{modicum} is simply $\EA$ extended by a function-symbol for each definable function in $\EA$. 
As usual, we are considering the usual metamathematical notions, in particular we use $\overline{n}$ to denote the (efficient \cite[p. 29]{buss1986bounded}) numeral of $n$; we also use $\#\varphi$ to denote the code of a formula, and $\ul\varphi\ur$ to denote its numeral. $\vec{Q}$ denotes an arbitrary array of quantifiers, $\vec{x}$ an array of variables; we define $\overline{\ex}:=\all$, $\overline{\all}:=\ex$, and $\overline{\vec{Q}}=\overline{Q_0\cdots Q_n}:=\overline{Q_0}\cdots\overline{Q_n}$. We define $|\vec{Q}|$ to be the length of the array, i.e. $|Q_0\cdots Q_n|=n+1$. For an array of quantifiers $\vec{Q}:=Q_0\cdots Q_n$ and an array of variables $\vec{x}:=x_0\cdots x_n$, clearly $\vec{Q}\vec{x}$ represents $Q_0x_0.\cdots Q_nx_n$. $Q^n$ denotes $\overbrace{Q\cdots Q}^{\text{$n$-times}}$. 

\begin{Definition} $\Formula(x)$ stands for a formula in $S$ that identifies all formulas, \emph{scilicet} such that $S\vdash \Formula(\overline{n})$ if, and only if, $n$ is the code of a formula. We also consider that $n\varr(x)$ defines all formulas with exactly $n$ free variables. As usual, $\all \Phi(x).\varphi(x)$ abbreviates $\all x. \Phi(x)\imp \varphi(x)$. In addition, $\l x_0,\ldots,x_k\r$ denotes a sequence function such that $S\vdash {(\l x_0,\ldots,x_k\r)}_i=x_i$. We assume that $S$ has function-symbols $\sub$ and $\numeral$ such that:
\begin{align*} &S\vdash \sub(\ul\varphi(x)\ur,\ul t\ur)=\ul\varphi(t)\ur;
&S\vdash \numeral(\overline{n})=\ul\overline{n}\ur.
\end{align*}
Feferman's dot notation (see \cite[p.\hspace{0.1cm}837]{Barwise1993} and \cite[p.\hspace{0.1cm}135]{buss1986bounded}) is defined in the following way: $\ul\varphi(\overset{\bigcdot}{x})\ur=\s(\ul\varphi\ur,x):=\sub(\ul\varphi\ur,\numeral(x))$. We consider function-symbols $\overbrace{\all x_0.}^{\bigcdot}$, $\overbrace{\ex x_0.}^{\bigcdot}$, $\overset{\bigcdot}{\neg}$, $\landd$, $\impd$, and $\eqd$ in $S$ satisfying:
\begin{align*} &S\vdash \overbrace{\all x_0.}^{\bigcdot}\ul\varphi(x_0)\ur=\ul\all x_0.\varphi(x_0)\ur; &&S\vdash \overbrace{\ex x_0.}^{\bigcdot}\ul\varphi(x_0)\ur=\ul\ex x_0.\varphi(x_0)\ur;\\
&S\vdash \overset{\bigcdot}{\neg}\ul\varphi\ur=\ul\neg\varphi\ur; &&S\vdash \ul\varphi\ur\landd\ul\psi\ur=\ul\varphi \land\psi\ur; \\
&S\vdash \ul\varphi\ur\lord\ul\psi\ur=\ul\varphi \lor\psi\ur; &&S\vdash \ul\varphi\ur\impd\ul\psi\ur=\ul\varphi \imp\psi\ur;\\
& S\vdash \ul\varphi\ur\eqd\ul\psi\ur=\ul\varphi \eq\psi\ur.
\end{align*}
We consider the previous definitions also for arrays of variables, namely $\ul\varphi(\overset{\bigcdot}{\vec{x}})\ur$ and $\overbrace{\vec{Q}\vec{x}.}^{\bigcdot}$. In fact, for an array of variables $T\vdash\sub(\ul\varphi(\vec{x})\ur,\ul\vec{t}\ur)=\ul\varphi(\vec{t})\ur$, where $\ul\vec{t}\ur:=\l\ul t_0\ur,\ldots, \ul t_n\ur\r$. Likewise, we also define $\s(y,\vec{x})=\s(y,x_0\cdots x_n):=\\\sub(y,\l \numeral(x_0),\ldots,\numeral(x_n)\r)$ (we are going to assume that $x_0\cdots x_n$ is the order of occurrences of the variables in $y$). If $\varphi$ is a sentence, we stipulate $T\vdash \sub(\ul\varphi\ur,\ul x_0\ur)=\ul\varphi\ur$; if $\varphi$ has more variables than the length $n$ of a given array $\vec{t}$, we define that $\sub$ applies to the first $n$ variables of $\varphi$ the array $\vec{t}$.
\end{Definition}
\begin{Definition} We say that \emph{$\alpha$ numerates $T$ \textup{(}in $S$\textup{)}}, or that \emph{$\alpha$ is a numeration of $T$}, if the set $\{\varphi|S\vdash\alpha(\ul\varphi\ur)\}$ coincides with the set of the codes of the axioms of $T$.
\end{Definition}

\begin{Definition} Given a numeration $\alpha$ of $T$, we consider the \emph{standard proof predicate for $\alpha$}, $\Prf_\alpha(x,y)$ that expresses ``$y$ is a proof of $x$ using axioms numerated by $\alpha$'', as defined by Feferman in \cite{Feferman1960} and we define the \emph{standard provability predicate for $\alpha$} by $\Pr_\alpha(x):=\ex y.\Prf_\alpha(x,y)$. Moreover, $\Proof_\alpha(x)$ expresses ``$x$ is a proof with axioms in $\alpha$'' (see \cite{Feferman1960}). As usual, $\Con_\alpha:=\neg \Pr_\alpha(\ul\perp\ur)$. Clearly, $S\vdash \Prf_\alpha(x,y)\imp \Proof_\alpha(y)$. We consider a fixed numeration $\alpha$ of $T$ in $S$ and we use the notation $\Pr_T$, $\Prf_T$, $\Con_T$, and so on.
\end{Definition}

We now present some of the derivability conditions satisfied by $\Pr_T$, for more details we recommend: \cite[pp.\hspace{0.1cm}133--149]{buss1986bounded}, \cite[pp.\hspace{0.1cm}117,\hspace{0.1cm}118]{buss1998handbook}, \cite{Feferman1960}, \cite[p. 14]{lind2017aspects}, \cite{Kurahashi_2018}, \cite{henk2017interpretability}, and \cite{kurahashi_2018aaa}. 

\begin{Fact}\label{fact1} The following conditions hold for $\Pr_T$:
\begin{description}
\item[\textbf{C1:}] If $T\vdash \varphi(\vec{x})$, then $S\vdash \Pr_T(\ulcorner \varphi(\overset{\bigcdot}{\vec{x}})\urcorner)$;
\item[\textbf{C2:}] $S\vdash \Pr_T(x\impd y)\imp(\Pr_T(x)\imp \Pr_T(y))$;
\item[\textbf{C3:}] For all ${\Sigma_1}(S)$-formulas $\varphi(\vec{x})$, $S\vdash \varphi(\vec{x})\imp \Pr_T(\ulcorner \varphi(\overset{\bigcdot}{\vec{x}})\urcorner)$;
\item[\textbf{C4:}] $S \vdash \alpha(x)\imp \Pr_T(x)$;
\item[\textbf{C6:}] $S\vdash \Pr_T(\ul\perp\ur\impd x)$;
\item[\textbf{C7:}] $S\vdash \ex x.\Pr_T(\s(y,x))\imp \Pr_T(\overbrace{\ex x_0.}^{\bigcdot}\sub(y,\ul x_0\ur))$;
\item[\textbf{C8:}] $S\vdash \Pr_T(\overbrace{\all x_0.}^{\bigcdot} \sub(y,\ul x_0\ur))\imp \all x. \Pr_T(\s(y,x))$.
\end{description}
\end{Fact}

\section{Our results}
We present a family of consistency statements that are similar to $\omega$-consistency; we do this to show that the actual array of quantifiers used in the $\omega$-consistency statement does not play a major role in some situations.

\begin{Definition} We define:
\begin{align*} &\Oca^{\vec{Q}}_T(y):=(\vec{Q}\vec{x}.\Pr_T(\s(y,\vec{x})))\imp\neg \Pr_T(\overbrace{\overline{\vec{Q}}\vec{x_0}.}^{\bigcdot}\negd \sub(y,\ul\vec{x_0}\ur));\\
&\Oca^{\vec{Q}}_T:=\all |\vec{Q}|\varr(y). \Oca^{\vec{Q}}_T(y).
\end{align*}
\end{Definition}

The next result confirms that $\Oca_T$ is a particular case of our generalized notion for different arrays of quantifers. 
\begin{Proposition}\label{prop5} $S\vdash \Oca_T\eq\Oca^{\all}_T$.
\end{Proposition}

\begin{proof} Immediate by definition. 
\end{proof}

The standard consistency statement is also a particular case of our analysis. 
\begin{Proposition}\label{prop6} $S\vdash \Oca^{\ex}_T\eq\Con_T$.
\end{Proposition}

\begin{proof} Reason inside $S$. Clearly, $\neg \Con_T\imp \neg \Oca^{\ex}_T$. Let us reason for the converse. Suppose that $\neg \Oca^{\ex}_T$. Then, $\ex 1\varr(y). \neg\Oca^{\ex}_\alpha(y)$. Take such a $y$. Consequently, $\ex x.\Pr_T(\s(y,{x}))\land \Pr_T(\overbrace{\all {x_0}.}^{\bigcdot}\negd \sub(y,\ul{x_0}\ur))$, and thus, using \textbf{C7}, we can conclude that $\Pr_T(\overbrace{\ex x_0.}^{\bigcdot}\sub(y,\ul x_0\ur))\land \\\Pr_T(\overbrace{\all {x_0}.}^{\bigcdot}\negd \sub(y,\ul{x_0}\ur))$, which entails $\Pr_T(\ul\perp\ur)$, i.e. $\neg \Con_T$.
\end{proof}

\begin{Definition} We stipulate that $\vec{Q_0}\sqsubseteq_0 \vec{Q_1}$ if there is an array $\vec{Q_2}$ such that $\vec{Q_1}=\vec{Q_2}\vec{Q_0}$ and $\vec{Q_0}\sqsubseteq_1 \vec{Q_1}$ if there is an array $\vec{Q_2}$ such that $\vec{Q_1}=\vec{Q_0}\vec{Q_2}$.
\end{Definition}

\begin{Theorem}\label{theo4} If $\vec{Q_0}\sqsubseteq_0 \vec{Q_1}$, then $S\vdash \Oca^{\vec{Q_1}}_T\imp \Oca^{\vec{Q_0}}_T$.
\end{Theorem}

\begin{proof} Reason inside $S$ and assume $\Oca^{\vec{Q_1}}_T$. Let $\vec{Q_2}$ be such that $\vec{Q_1}=\vec{Q_2}\vec{Q_0}$. Consider $y$ such that $|\vec{Q_0}|\varr(y)$ and $\vec{Q_0}\vec{x}.\Pr_T(\s(y,\vec{x}))$. Take
\begin{align*} y':=y\landd \ul\small\bigwedge_{i=0}^{|\vec{Q_2}|-1} v_i=v_i\ur,
\end{align*}
where each $v_i$ does not occur in $y$. Then, considering $\vec{v}:=v_0\cdots v_{|\vec{Q_2}|-1}$ and $\vec{x'}:=\vec{x}\vec{v}$, we can conclude $\vec{Q_2}\vec{v}.\vec{Q_0}\vec{x}.\Pr_T(\s(y',\vec{x'}))$. From the assumption of $\Oca^{\vec{Q_1}}_T$, it follows $\neg \Pr_T(\overbrace{\overline{\vec{Q_2}}\vec{v}.\overline{\vec{Q_0}}\vec{x}.}^{\bigcdot}\negd \sub(y',\ul\vec{x'}\ur))$. As
\begin{align*} \Pr_T\biggl(\hspace{-0.1cm}(\overbrace{\overline{\vec{Q_2}}\vec{v}.\overline{\vec{Q_0}}\vec{x}.}^{\bigcdot}\negd \sub(y',\ul\vec{x'}\ur))\hspace{-0.1cm}&\eqd \hspace{-0.1cm}\Big(\overbrace{\overline{\vec{Q_2}}\vec{v}.\overline{\vec{Q_0}}\vec{x}.}^{\bigcdot}\negd\sub(y,\ul\vec{x}\ur)\lord \ul\hspace{-0.1cm}\small\bigvee_{i=0}^{|\vec{Q_2}|-1} \hspace{-0.1cm}\neg v_i=v_i\ur\hspace{-0.08cm}\Big)\\
&\eqd\overbrace{\overline{\vec{Q_0}}\vec{x}.}^{\bigcdot}\negd \sub(y,\ul\vec{x}\ur)\biggl)
\end{align*}
we can therefore conclude $\neg \Pr_T(\overbrace{\overline{\vec{Q_0}}\vec{x}.}^{\bigcdot}\negd \sub(y,\ul\vec{x}\ur))$.
\end{proof}

\begin{Theorem}\label{theo5} If $\vec{Q_0}\sqsubseteq_1 \vec{Q_1}$, then $S\vdash \Oca^{\vec{Q_1}}_T\imp \Oca^{\vec{Q_0}}_T$.
\end{Theorem}

\begin{proof} Similar to the previous proof. 
\end{proof}

\begin{Proposition}\label{Proposition8} For any array of quantifiers $\vec{Q}$, $S\vdash \Oca^{\all\all\vec{Q}}_T\eq \Oca^{\all\vec{Q}}_T$.
\end{Proposition}

\begin{proof} Reason inside $S$. From the previous result, $\Oca^{\all\all\vec{Q}}_T\imp \Oca^{\all\vec{Q}}_T$. Let us argue for the converse. Assume $\Oca^{\all\vec{Q}}_T$ and that $y$ is such that it satisfies $|\all\all\vec{Q}|\varr(y)$ and 
\begin{align*}\tag{I}\all x_0.\all x_1.\vec{Q}\vec{x}.\Pr_T(\s(y,x_0 x_1\vec{x})). 
\end{align*}
Now, consider $y'$ as being $\sub(y,\ul {(v)}_0{(v)}_1\vec{x}\ur)$; this means that $y'$ has variables $v$ and all the $\vec{x}$, and that the occurrences of $x_0$ in $y$ are occurrences of ${(v)}_0$ in $y'$, likewise, the occurrences of $x_1$ are substituted for occurrences of ${(v)}_1$. 

Take any arbitrary $v$. From (I), $\vec{Q}\vec{x}.\Pr_T(\s(y,{(v)}_0 {(v)}_1\vec{x}))$, consequently\footnote{Because $S\vdash \Pr_T(\ul\overbrace{{(v)}_i}^{\bigcdot}={(\overset{\bigcdot}{v})}_i\ur)$, with $i\in\{0,1\}$.} we obtain\\ $\all v.\vec{Q}\vec{x}.\Pr_T(\s(y',v\vec{x}))$. From $\Oca^{\all\vec{Q}}_T$, $\neg \Pr_T(\overbrace{\all v.\overline{\vec{Q}}\vec{x}.}^{\bigcdot}\negd \sub(y',\ul v\vec{x}\ur))$. As\footnote{It is a routine matter to prove this; the idea is, when reasoning inside $S\vdash \Pr_T(\cdot)$, that for the left to right implication we consider $v=\l x_0,x_1\r$, and for the converse direction we take $x_0={(v)}_0$ and $x_1={(v)}_1$.}
\begin{align*} \Pr_T((\overbrace{\all v.\overline{\vec{Q}}\vec{x}.}^{\bigcdot}\negd \sub(y',\ul v\vec{x}\ur))\eqd(\overbrace{\all x_0.\all x_1\overline{\vec{Q}}\vec{x}.}^{\bigcdot}\negd \sub(y,\ul x_0x_1\vec{x}\ur))),
\end{align*}
then
\begin{align*} \Pr_T(\overbrace{\all v.\overline{\vec{Q}}\vec{x}.}^{\bigcdot}\negd \sub(y',\ul v\vec{x}\ur))\eq\Pr_T(\overbrace{\all x_0.\all x_1\overline{\vec{Q}}\vec{x}.}^{\bigcdot}\negd \sub(y,\ul x_0x_1\vec{x}\ur)),
\end{align*}
and we conclude $\neg \Pr_T(\overbrace{\all x_0.\all x_1\overline{\vec{Q}}\vec{x}.}^{\bigcdot}\negd \sub(y,\ul x_0x_1\vec{x}\ur))$, as wanted.
\end{proof}

\begin{Proposition}\label{proposi9} For any array of quantifiers $\vec{Q}$, $S\vdash \Oca^{\ex\ex\vec{Q}}_T\eq \Oca^{\ex\vec{Q}}_T$.
\end{Proposition}

\begin{proof} Similar to the previous proof.
\end{proof}

\begin{Proposition}\label{props10} $S\vdash \Oca^{\all}_T\eq \Oca^{\all\ex}_T$.
\end{Proposition}

\begin{proof} It is easy to conclude that $S\vdash \Oca^{\all\ex}_T\imp \Oca^{\all}_T$ using the reasoning of Theorem \ref{theo4}. Reason inside $S$ for the converse. Assume $\Oca^{\all}_T$ and that $y$ is such that it satisfies $|\all\ex|\varr(y)$ and $\all x_0.\ex x_1.\Pr_T(\s(y,x_0x_1))$. We may assume, without loss of generality, that $x_0$ and $x_1$ are exactly the variables of $y$ (this avoids occurrences of `$\sub$', since $S\vdash y=\sub(y,\ul x_0x_1\ur)$). Then, using \textbf{C7}, $\all x_0.\Pr_T(\s(\overbrace{\ex x_1.}^{\bigcdot}y,x_0))$. From $\Oca^{\all}_T$, $\neg \Pr_T(\overbrace{\ex x_0.}^{\bigcdot}\negd \overbrace{\ex x_1.}^{\bigcdot}y)$, and so $\neg \Pr_T(\overbrace{\ex x_0.}^{\bigcdot} \overbrace{\all x_1.}^{\bigcdot}\negd y)$; as wanted.
\end{proof}~\\
\vspace{-0.9cm}
\begin{Proposition}\label{propost11} $S\vdash \Oca^{\all}_T\eq \Oca^{\ex\all}_T$.
\end{Proposition}

\begin{proof} From Theorem \ref{theo4}, $S\vdash \Oca^{\ex\all}_T\imp \Oca^{\all}_T$. Reason inside $S$ for the converse. Assume $\Oca^{\all}_T$ and that $y$ is such that it satisfies $|\ex\all|\varr(y)$ and $\ex x_0.\all x_1.\Pr_T(\s(y,x_0x_1))$. Similarly to what was previously done, we may assume, without loss of generality, that $x_0$ and $x_1$ are exactly the variables of $y$. Consider $y'$ as being $y$ where the occurrence of $x_0$ is substituted by the numeral of $x_0$, i.e. $y'=\s(y,x_0)$, assuming $x_0$ is the first variable occurring in $y$. Then, $\ex x_0.\all x_1.\Pr_T(\s(y',x_1))$. Clearly, $\all x_0.1\varr(y')$, and so, by $\Oca^{\all}_T$, $\all x_0.((\all x_1.\Pr_T(\s(y',x_1)))\imp \neg \Pr_T(\overbrace{\ex x_1.}^{\bigcdot}\negd y'))$. Thus, $\ex x_0. \neg\Pr_T(\overbrace{\ex x_1.}^{\bigcdot}\negd y')$. As we know from \textbf{C8},\\ $\Pr_T(\overbrace{\all x_0.}^{\bigcdot}\overbrace{\ex x_1.}^{\bigcdot}\negd y)\imp \all x_0.\Pr_T(\overbrace{\ex x_1.}^{\bigcdot}\negd y')$, we get $\neg \Pr_T(\overbrace{\all x_0.}^{\bigcdot}\overbrace{\ex x_1.}^{\bigcdot}\negd y)$, as wanted.
\end{proof}

\begin{Proposition} For any $\vec{Q}$, $S\vdash \Oca^{\vec{Q}\ex}_T\eq \Oca^{\vec{Q}}_T$.
\end{Proposition}

\begin{proof} Reason inside $S$. From Theorem \ref{theo5}, $\Oca^{\vec{Q}\ex}_T\imp \Oca^{\vec{Q}}_T$; so, let us argue for the other implication. Suppose $\Oca^{\vec{Q}}_T$. Furthermore, assume the antecedent of $\Oca^{\vec{Q}\ex}_T$, that is to say, consider $y$ such that $|\vec{Q}\ex|\varr(y)$ and $\vec{Q}\vec{x}.\ex z.\Pr_T(\s(y,z,\vec{x}))$. Assume, without loss of generality, that $z$ is the first variable occurring in the formula $y$ (this assumption is made to simplify the way we are considering the substitution function). From the derivability conditions, $\vec{Q}\vec{x}.\Pr_T(\overbrace{\ex z.}^{\bigcdot}\s(\sub(y,\ul z\ur),\vec{x}))$, and thus $\vec{Q}\vec{x}.\Pr_T(\s(\overbrace{\ex z.}^{\bigcdot}\sub(y,\ul z\ur),\vec{x}))$. Clearly, $|\vec{Q}|\varr(\overbrace{\ex z.}^{\bigcdot}\sub(y,\ul z\ur))$; consequently, from the assumption of $\Oca^{\vec{Q}}_T$, we can conclude that 
\begin{align*}\neg \Pr_T(\overbrace{\overline{\vec{Q}}\vec{x_0}.}^{\bigcdot}\negd \sub(\overbrace{\ex z.}^{\bigcdot}\sub(y,\ul z\ur),\ul\vec{x_0}\ur)),
\end{align*}
and so
\begin{align*}\neg \Pr_T(\overbrace{\overline{\vec{Q}}\vec{x_0}.\overline{\ex} z.}^{\bigcdot}\negd \sub(y,\ul z\ur,\ul\vec{x_0}\ur)),
\end{align*}
as desired; this confirms $\Oca^{\vec{Q}\ex}_T$.
\end{proof}

\begin{Main}\label{main} For every $n,m\in\mathbb{N}$, $S\vdash\Oca_T\eq \Oca^{\all^{n+1}}_T\eq\Oca^{\all^{n+1}\ex^m}_T\eq\Oca^{\ex\all\ex^n}_T$.
\end{Main}

\begin{proof} Immediate from the preceding propositions.
\end{proof}

Now we move to create a theory of truth that encompasses every statement $\Oca^{\vec{Q}}_T$. \newpage

\begin{Definition} For each $T$, we define the theory $\Tr(S)_T$ has being the theory obtained from $S$ extended by a predicate $\Tr(x)$ and the following axioms:
\begin{description}
\item[A1] $\Pr_T(x)\imp \Tr(x)$;
\item[A2] $\Tr(\overbrace{\all x_0.}^{\bigcdot} \sub(y,\ul x_0\ur))\eq \all x. \Tr(\s(y,x))$;
\item[A3] $\Tr(\overbrace{\ex x_0.}^{\bigcdot} \sub(y,\ul x_0\ur))\eq \ex x. \Tr(\s(y,x))$;
\item[A4] $\Tr(\negd y)\eq \neg\Tr(y)$;
\item[A5] $\Tr(x)\imp \Formula(x)$.
\end{description}
\end{Definition}

\begin{Lemma} $\Tr(S)_T\vdash \Pr_T(\overbrace{\vec{Q}\vec{x_0}.}^{\bigcdot}\sub(y,\ul\vec{x_0}\ur))\imp \vec{Q}\vec{x}.\Tr(\s(y,\vec{x}))$.
\end{Lemma}

\begin{proof} Let $\vec{Q}:=Q_0\cdots Q_n$ and $\vec{x}:=x_0 \cdots x_n$. Reason inside $\Tr(S)_T$. Assume that\\ $\Pr_T(\overbrace{\vec{Q}\vec{x}.}^{\bigcdot}\sub(y,\ul\vec{x}\ur))$. We may assume, without loss of generality, that $\vec{x}$ are the variables in $y$. So, $\Pr_T(\overbrace{\vec{Q}\vec{x}.}^{\bigcdot}y)$. By \textbf{A1}, we can conclude that $\Tr(\overbrace{Q_0x_0.\cdots Q_nx_n.}^{\bigcdot}y)$. By several applications of \textbf{A2} and \textbf{A3} we conclude
\begin{align*} \Tr(\overbrace{Q_0x_0.\cdots Q_nx_n.}^{\bigcdot}y) &\imp\Tr(\overbrace{Q_0x_0.}^{\bigcdot}\sub(\overbrace{Q_1x_1.\cdots Q_nx_n.}^{\bigcdot}y,\ul x_0\ur)\\
&\imp Q_0x_0. \Tr(\s(\overbrace{Q_1x_1.\cdots Q_n x_n.}^{\bigcdot}y,x_0))\imp \cdots\\
&\cdots \imp Q_0x_0.\cdots Q_nx_n.\Tr(\s(y,x_0\cdots x_n)),
\end{align*}
as wanted.
\end{proof}

\begin{Theorem}\label{theo55} For any array of quantifiers $\vec{Q}$, $\Tr(S)_T\vdash \Oca^{\vec{Q}}_T$.
\end{Theorem}

\begin{proof} Reason inside $\Tr(S)_T$. Let $y$ be such that $|\vec{Q}|\varr(y)$ and also \newline$\vec{Q}\vec{x}.\Pr_T(\s(y,\vec{x}))$. Suppose, aiming at establishing a contradiction, that we have\\ $\Pr_T(\overbrace{\overline{\vec{Q}}\vec{x_0}.}^{\bigcdot}\negd \sub(y,\ul\vec{x_0}\ur))$, i.e. $\Pr_T(\overbrace{\overline{\vec{Q}}\vec{x_0}.}^{\bigcdot} \sub(\negd y,\ul\vec{x_0}\ur))$. By \textbf{A1}, from the initial assumption we get $\vec{Q}\vec{x}.\Tr(\s(y,\vec{x}))$. By the previous Lemma applied to the second assumption, we obtain $\overline{\vec{Q}}\vec{x}.\Tr(\s(\negd y,\vec{x}))$. Thus, using \textbf{A4}, we can conclude $\vec{Q}\vec{x}.\Tr(\s(y,\vec{x})). \land \overline{\vec{Q}}\vec{x}.\neg \Tr(\s(y,\vec{x}))$, and so $\vec{Q}\vec{x}.\Tr(\s(y,\vec{x})).\land \neg \vec{Q}\vec{x}.\Tr(\s(y,\vec{x}))$, a contradiction.
\end{proof}

\section{Conclusions}  
We introduced a generalization of the usual way of internalizing the $\omega$-consistency statement by allowing the occurrence of different arrays of quantifiers $\vec{Q}$; we denoted it by $\Oca^{\vec{Q}}_T$. We showed that for some specific arrays of quantifiers we could, provably in $S$, retrieve the original $\omega$-consistency; moreover, we developed a theory of truth that proves each statement $\Oca^{\vec{Q}}_T$.

The next figure summarizes the main results of our paper, where we introduced new kinds of $\omega$-consistency.
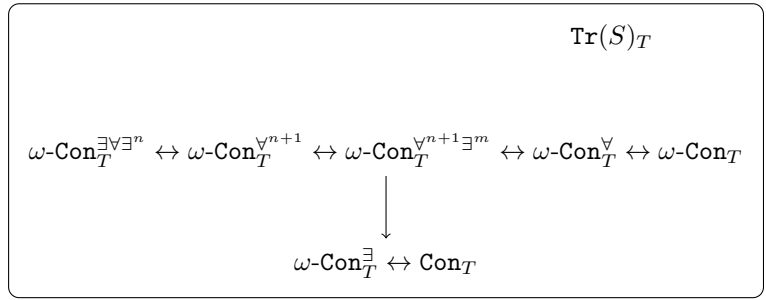
\begin{figure}[H]
\centering
\begin{tikzpicture}[framed,grow = up, edge from parent/.style={draw, <-,},sibling distance=60pt,background rectangle/.style={draw=black, rounded corners}]
\node{$\Oca^{\ex}_T\eq \Con_T$}
	child {node{$\Oca^{\ex\all\ex^n}_T\eq\Oca^{\all^{n+1}}_T\eq\Oca^{\all^{n+1}\ex^m}_T\eq\Oca^{\all}_T\eq \Oca_T$}};
\node at (3,3) {$\Tr(S)_T$};
\end{tikzpicture}
\caption{$S$-implications and $S$-equivalences for the general notion of consistency $\Oca^{\vec{Q}}_T$ depending on the array of quantifiers; $S$-implications are read from the top to the bottom, \emph{id est} the upper part implies the lower part; all the formulas are provable in $\Tr(S)_T$. This image relies on our Main Theorem and Theorem \ref{theo55}.}
\end{figure}

\section{Acknowledgement}
The author was supported by national funds through the FCT\textendash Funda\c{c}\~{a}o para a Ci\^{e}ncia e a Tecnologia, I.P., under the scope of the projects UIDB/04561/2020 and UIDP/04561/2020 (CMAFcIO, UL).

\bibliographystyle{abbrv}
\bibliography{cons}
\end{document}